\documentclass[11pt]{article} 
 
\usepackage{latexsym} 
\usepackage{amsfonts} 
\usepackage{amsmath} 
\usepackage{amsthm} 
\usepackage{amssymb}
\usepackage{mathrsfs} 
\usepackage{dsfont} 
\usepackage{bbold} 
\usepackage[english]{babel} 
\usepackage{caption} 
\usepackage{epsfig} 
\usepackage{float} 
\usepackage{subfigure} 
\usepackage{psfrag} 
\usepackage{graphicx} 
\usepackage{epsfig} 
\usepackage{amsbsy} 
\usepackage{hyperref}
 
\DeclareMathOperator{\Val}{\matV} 
 
 \DeclareMathOperator{\diag}{diag} 

\newtheorem{theorem}{Theorem} 
\newtheorem*{prop*}{Theorem} 
\newtheorem{theo}{Theorem} 
 
\newtheorem{defi}[theorem]{Definition} 
\newtheorem{lemma}[theorem]{Lemma} 
 
\newtheorem{rmk}[theorem]{Remark}

\newtheorem{hyp}{Hypothesis}

\newcommand{\zerarcounters}{\setcounter{equation}{0}\setcounter{theorem}{0}} 
 
\newcommand{\beq}{\begin{equation}}
\newcommand{\eeq}{\end{equation}}

\newcommand{\ZZZ}{\mathds{Z}} 
\newcommand{\CCC}{\mathds{C}} 
\newcommand{\NNN}{\mathds{N}} 
 
\newcommand{\RRR}{\mathds{R}} 
\newcommand{\TTT}{\mathds{T}} 
\newcommand{\uno}{\mathds{1}} 
 
\newcommand{\CCCC}{{\mathcal C}}

\newcommand{\calF}{{\mathcal F}} 
\newcommand{\calG}{{\mathcal G}}

\newcommand{\MM}{{\mathcal M}} 
\newcommand{\NN}{{\mathcal N}} 
 
\newcommand{\calP}{{\mathcal P}}

\newcommand{\gotB}{{\mathfrak B}}

\newcommand{\gotN}{{\mathfrak N}}

\newcommand{\gotR}{{\mathfrak R}}

\newcommand{\matV}{{\mathscr V}}

\newcommand{\ol}{\overline} 
 
\newcommand{\Fullbox}{{\rule{2.0mm}{2.0mm}}} 
 
\newcommand{\EP}{\hfill\Fullbox\vspace{0.2cm}} 
\newcommand{\prova}{\noindent{\it Proof. }} 
\newcommand{\io}{\infty} 
\newcommand{\e}{\varepsilon} 
\newcommand{\al}{\alpha} 
\newcommand{\de}{\delta} 
\newcommand{\be}{\beta}

\newcommand{\x}{\xi} 
\newcommand{\p}{\pi}

\newcommand{\om}{\omega}

\newcommand{\la}{\lambda} 
\newcommand{\f}{\varphi}

\newcommand{\del}{\partial}

\newcommand{\oo}{{\omega}}

\newcommand{\nn}{{\nu}}

\newcommand{\vzero}{{0}}

\newcommand{\ii}{{\rm i}}

\def\ins#1#2#3{\vbox to0pt{\kern-#2 \hbox{\kern#1 #3}\vss}\nointerlineskip} 
 
\begin{document}
 
%%%%%%%%%%%%%%%%%%%%%%%%%%%%%%%%%%%%%%%%%%%%%%%%%
%\title{\bf Response solutions of arbitrary codimension}
\title{\bf Resonant tori of arbitrary codimension for quasi-periodically forced systems}
%%%%%%%%%%%%%%%%%%%%%%%%%%%%%%%%%%%%%%%%%%%%%%%%%
 
\author{\bf Livia Corsi$^{1}$\footnote{Supported by the European Research Council under FP7
``Hamiltonian PDEs and small divisor problems: a dynamical systems approach".}
\ and Guido Gentile$^{2}$
\vspace{2mm} 
\\ \small
$^{1}$ Dipartimento di Matematica, Universit\`a di Roma ``La Sapienza", Roma, 00185, Italy
\\ \small 
$^{2}$ Dipartimento di Matematica e Fisica, Universit\`a di Roma Tre, Roma, 00146, Italy
\\ \small 
E-mail:  corsi@mat.uniroma1.it, gentile@mat.uniroma3.it}

\date{} 
 
\maketitle 

%%%%%%%%%%%%%%%%%%%%%%%%%%%%%%%%%%%%%%%%%%%%%%%%%
\begin{abstract} 
We consider a system of rotators subject to a small quasi-periodic forcing.
We require the forcing to be analytic and satisfy a time-reversibility property
and we assume its frequency vector to be Bryuno. Then we prove that, without imposing
any non-degeneracy condition on the forcing, there exists at least one quasi-periodic
solution with the same frequency vector as the forcing. The result can be interpreted as a theorem
of persistence of lower-dimensional tori of arbitrary codimension in degenerate cases.
\end{abstract} 
%%%%%%%%%%%%%%%%%%%%%%%%%%%%%%%%%%%%%%%%%%%%%%%%%
  
%\tableofcontents

%%%%%%%%%%%%%%%%%%%%%%%%%%%%%%%%%%%%%%%%%%%%%%%%%
%%%%%%%%%%%%%%%%%%%%%%%%%%%%%%%%%%%%%%%%%%%%%%%%%
\zerarcounters 
\section{Introduction} 
\label{sec:1} 
%%%%%%%%%%%%%%%%%%%%%%%%%%%%%%%%%%%%%%%%%%%%%%%%%
%%%%%%%%%%%%%%%%%%%%%%%%%%%%%%%%%%%%%%%%%%%%%%%%%

In this paper we deal with the problem of persistence of lower-dimensional tori for
nearly integrable Hamiltonian systems without assumptions of non-degeneracy
on the perturbation. Consider an analytic Hamiltonian function of the form
\begin{equation} \label{eq:1.1}
H(\f,J)= H_{0}(J) + \e f(\f,J) ,
\end{equation}
where $(J,\varphi) \in \RRR^{n} \times \TTT^{n}$ are conjugate action-angle variables,
with $\TTT=\RRR/2\p\ZZZ$, and $\e\in\RRR$. The function $H_0$ will be referred to
as the unperturbed Hamiltonian, while $f$ will be called the perturbation.
For $\e=0$ the phase space is foliated into $n$-dimensional invariant tori;
on each torus the action $J$ is constant and the motion is quasi-periodic
with frequency vector $\oo_{0}(J):=\del_{J}H_{0}(J)$.

We say that a frequency vector $\oo_0\in\RRR^{n}$ is resonant with multiplicity $r$
if there is a rank $r$ subgroup $G$ of $\ZZZ^{n}$ such that $\oo_0\cdot\nn=0$ for all $\nn\in G$
and $\oo_0\cdot\nn \neq 0$ for all $\nn \in \ZZZ^{n}\setminus G$.
It is a long standing conjecture that, in nearly integrable Hamiltonian systems
with Hamiltonian function of the form \eqref{eq:1.1}, under a convexity assumption on $H_0$,
for any $r\le n-1$ and for most families of unperturbed tori whose frequency vectors
are resonant with multiplicity $r$, at least $r+1$ tori survive any perturbation for
$\e$ small enough \cite{CKLY,KP}. When this happens one says that the surviving torus
has codimension $r$. In the case $r=n-1$, where the surviving tori reduce to
closed orbits, this has been shown by Bernstein and Katok \cite{BK}.
The case $r<n-1$ is harder, because of the presence of
small divisors. For $r=1$ the conjecture has been proved by Cheng \cite{C2},
while only partial results exist for $1<r<n-1$, requiring non-degeneracy hypotheses
on the perturbation \cite{CW,T,G1}.

A somewhat different problem is the following. Let the frequency vector $\oo_0$ be resonant
with multiplicity $r$ and fix $J_0$ such that $\oo_{0}(J_0)=\oo_0$.
Then, after a suitable change of coordinates, one can rewrite \eqref{eq:1.1} as
\begin{equation} \label{eq:1.2}
H(\al,\be,A,B)= H_{0}(A,B) + \e f(\al,\be,A,B) ,
\end{equation}
where $(A,B)\in \RRR^{d}\times\RRR^{r}$ and $(\al,\be) \in \TTT^{d}\times\TTT^{r}$
are such that $d+r=n$ and $\oo_0(A_0,B_0)=(\oo,0)$, with $\oo\in\RRR^{d}$ non-resonant;
not to overwhelm the notation, we are denoting with the same symbol the Hamiltonian
obtained by composing the original Hamiltonian with the change of coordinates.
Assume some Diophantine condition on $\oo$ and consider the particular family of tori
of the unperturbed system with that frequency vector. Then look at the fate of the family
when the perturbation is switched on. For Hamiltonian functions \eqref{eq:1.2},
again assuming $H_0$ to be convex, in general one expects only that at least one torus
in the family survives any perturbation $f$ for $\e$ small enough. In fact, 
the persistence of one such torus is known to hold for any $r$ only
if one assumes a non-degeneracy condition on $f$. More precisely, the perturbation averaged
on the unperturbed torus is taken to have a non-degenerate maximum \cite{T};
see also \cite{JDZ,GG}, where the analyticity
properties of the tori in the perturbation parameter are investigated.
By weakening the non-degeneracy assumptions on the perturbation
partial results have been obtained for $r=1$ in \cite{Y,GGG}.

Despite the fact that assumptions on the perturbations are commonly believed to be unnecessary, 
so far the only result existing in the literature  for Hamiltonian functions \eqref{eq:1.2}
is in the case of tori with codimension 1 and is due to Cheng \cite{C1}. 
An analogous result of persistence of at least one lower-dimensional torus
of codimension 1, without assuming any hypothesis on the perturbation, has been recently proved
for a suitable class of non-convex $H_0$ with a saddle point \cite{PK,CFG}.

In the present paper we consider a system with ``partially isochronous" unperturbed Hamiltonian
and a perturbation depending only on the angles, namely
\begin{equation} \label{eq:1.3}
H(\al,\be,A,B)= \om \cdot A + \frac{1}{2} B^{2} + \e f(\al,\be) ,
\end{equation}
and, under a mild Diophantine condition on $\om\in\RRR^{d}$,
we prove the existence of at least one quasi-periodic solution with frequency vector
$\om$ for the corresponding Hamiltonian system, requiring on $f$ a parity
condition with respect to the variable $\al$; see Hypothesis \ref{hyp:2} below.
Such assumption can be seen as a time-reversibility condition on the Hamiltonian system.

For non-convex Hamiltonians \eqref{eq:1.3} with $r=1$, the same result of persistence
of at least one invariant torus has been proved in \cite{CG1} without assuming
Hypothesis \ref{hyp:2} on the perturbation; the result has been extended to more
general one-dimensional systems in \cite{CG2}.
Strictly speaking, Cheng's result does not imply the result in \cite{CG1,CG2},
because the unperturbed Hamiltonian is not convex; however Cheng's method can be
adapted to such a case; see \cite{F} for an explicit implementation.
Moreover in \cite{CG1,CG2} the frequency vector is allowed to satisfy the Bryuno condition,
a weaker condition with respect to the standard Diophantine condition considered
in both \cite{C1} and \cite{F}; in \cite{CG2} non-Hamiltonian systems have been
considered as well: in that case further conditions have to be assumed in general
on the perturbation (in the Hamiltonian case such conditions are automatically satisfied).

For the general case \eqref{eq:1.1} with $r>1$,
the problem of persistence of at least one lower-dimensional torus without requiring any
non-degeneracy assumption on the perturbation is still open. 
As far as we know, Theorem \ref{thm:1} below is the first result in that direction.
The time reversibility hypothesis we assume on the perturbation is a symmetry property,
not a non-degeneracy  condition, and aims to ensure a suitable cancellation that we need in the proof.
It is natural to conjecture such a cancellation to hold in general, but this seems to require further work.

%%%%%%%%%%%%%%%%%%%%%%%%%%%%%%%%%%%%%%%%%%%%%%%%%
\subsection{Statement of the results}
%%%%%%%%%%%%%%%%%%%%%%%%%%%%%%%%%%%%%%%%%%%%%%%%%

Let us write the Hamiltonian \eqref{eq:1.3} in the non-autonomous form
\begin{equation} \label{eq:1.4}
H(\be,B,t) = \frac{1}{2} B^2 + \e f(\oo t,\be) ,
\end{equation}
where $(B,\be)\in \RRR^{r}\times\TTT^{r}$ are conjugate action-angle variables,
$\oo\in\RRR^d$ is the \emph{frequency vector},
$f$ is analytic in the complexified torus $\TTT^{d+r}_\x:=\{\f\in\CCC^{d+r} : {\rm Re}\,\f_{k}
\in\TTT,\,|{\rm Im}\f_{k}|\le \x\}$ and $\e\in\RRR$ is a small parameter ({\it perturbation parameter}).
Without loss of generality, we may and shall assume that
$\oo$ has rationally independent components and $\e \ge 0$.

The corresponding Hamilton equations can be written as
\begin{equation} \label{eq:1.5}
\ddot{\be} =-\e \del_\be f(\oo t,\be) 
\end{equation}
and hence describe an $r$-dimensional system of rotators with a quasi-periodic forcing.

We consider $d>1$ in the following, since the periodic case $d=1$ is well understood;
see \cite{CG0} and references within. For $\om\in\RRR^{d}$ define
\begin{equation} \label{eq:1.6}
\alpha_{m}(\om) := \inf_{0<|\nu|\le 2^{m}} |\om\cdot\nu| , \qquad
\gotB(\om):=\sum_{m=0}^{\io} \frac{1}{2^{m}} \log \frac{1}{\alpha_{m}(\om)} ,
\end{equation}
where $\cdot$ denotes the inner product in $\RRR^{d}$ and
$|\nu|:=|\nu_{1}|+\ldots+|\nn_{d}|$, if $\nu_1,\ldots,\nu_d$ are the components of $\nu$.
We assume the following hypotheses on $\omega$ and $f$.

%%%%%%%%%%%%%%%%%%%%%%%%%%%%%%%%%%%%%%%%%%%%%%%%%
\begin{hyp} \label{hyp:1} 
$\om\in\RRR^{d}$, with $d>1$, satisfies the condition $\gotB(\om)<\io$.
\end{hyp}
%%%%%%%%%%%%%%%%%%%%%%%%%%%%%%%%%%%%%%%%%%%%%%%%%

%%%%%%%%%%%%%%%%%%%%%%%%%%%%%%%%%%%%%%%%%%%%%%%%%
\begin{hyp} \label{hyp:2} 
$f$ is even in $\al$, that is $f(-\al,\be)=f(\al,\be)$.
\end{hyp}
%%%%%%%%%%%%%%%%%%%%%%%%%%%%%%%%%%%%%%%%%%%%%%%%%

Hypothesis \ref{hyp:1} is a weak Diophantine condition on the frequency vector,
known in the literature as the \emph{Bryuno condition}. Hypothesis 2 is a
\emph{time reversibility} condition on the forcing.
We shall prove the following result.

%%%%%%%%%%%%%%%%%%%%%%%%%%%%%%%%%%%%%%%%%%%%%%%%%
\begin{theo} \label{thm:1}
Consider the Hamiltonian \eqref{eq:1.4} and assume Hypotheses \ref{hyp:1} and \ref{hyp:2}.
Then for $\e$ small enough there exists at least one quasi-periodic solution to the
corresponding Hamilton equation \eqref{eq:1.5} with frequency vector $\om$.
\end{theo}
%%%%%%%%%%%%%%%%%%%%%%%%%%%%%%%%%%%%%%%%%%%%%%%%%

%%%%%%%%%%%%%%%%%%%%%%%%%%%%%%%%%%%%%%%%%%%%%%%%%
\subsection{Informal presentation of the proof}\label{sec:1.3}
%%%%%%%%%%%%%%%%%%%%%%%%%%%%%%%%%%%%%%%%%%%%%%%%%

We look for quasi-periodic solutions to \eqref{eq:1.5} with frequency vector $\om$
({\it response solutions}),
so we split
\begin{equation} \nonumber
\be(t)=\be_0+b(\oo t), \qquad [b(\cdot)]_{0} := \frac{1}{(2\pi)^{d}} 
\int_{\TTT^{d}} {\rm d}\alpha \, b(\alpha) = 0 ,
\end{equation}
and separate \eqref{eq:1.5} into the so-called ``range'' and ``bifurcation'' equations
\begin{subequations}\label{eq:1.7}
\begin{align}
& (\om\cdot\del_\al)^2 b + \e \left(
\del_\be f(\al,\be_0+b) - [ \del_\be f(\cdot,\be_0+b(\cdot))]_{0} \right) =0,
\label{eq:1.7a} \\
& [\e \del_\be f(\cdot,\be_{0}+b(\cdot))]_{0}=0 .
\label{eq:1.7b}
\end{align}
\end{subequations}

The range equation \eqref{eq:1.7a} can be seen as an implicit function equation of the form
\begin{equation} \label{eq:1.8}
D_{\om} b + \e F(\be_0+b)=0,
\end{equation}
where the differential operator $D_{\om} :=(\om\cdot\del_\al)^2$  is diagonal on the
Fourier basis $\{{\rm e}^{\ii\nu\cdot\al}\}_{\nu\in\ZZZ^d}$. As it is well known, 
the standard implicit function theorem fails to apply since the inverse of $D_{\om}$ is unbounded 
(small divisor problem). In order to overcome this difficulty, one can implement a fast iterative scheme
(KAM, Nash-Moser or Renormalisation Group). At any iterative step
one has to impose Diophantine conditions involving the parameter $\be_0$, so as to control
the small divisors: at the $n$-step one requires $\be_0$ to be restricted to some set $\CCCC_{n} \subset \TTT^{r}$.
Eventually one obtains a function $b=b(\al,\e,\be_0)$ which solves the range equation
for $\be_0$ in a (possibly empty) set $\CCCC=\cap_{n=0}^{\io} \CCCC_{n}$.
Then one passes to the bifurcation equation \eqref{eq:1.7b}.
Formally (i.e. by assuming that a solution to \eqref{eq:1.7a} exists), by relying on
the variational structure of the Hamilton equation one has the identity
\begin{equation}\label{eq:1.9}
[\e\del_\be f(\cdot,\be_0+b(\cdot,\e,\be_0))]_0 = -\del_{\be_0} L(\e,\be_0),
\end{equation}
where $L(\e,\be_0)$ is the average with respect to $\al$ of the Lagrangian computed along
$b(\al,\e,\be_0)$; see for instance \cite{R,ACE,B}. However, without requiring any assumption
on the perturbation, it is difficult to have any control on the set $\CCCC$. Even assuming
that the solution can be extended to the torus $\TTT^{r}$, so as to ensure $L(\e,\be_0)$
to have a critical point, it is by no way obvious that such a critical point falls inside $\CCCC$.
In fact, the non-degeneracy condition usually assumed in the literature aims exactly to
make the left hand side of \eqref{eq:1.9} to vanish for a suitable value of  $\be_0$ for which
the Diophantine conditions are satisfied.
So, if no non-degeneracy condition is assumed on the perturbation, we have to proceed in a different way.

In the present paper, we use a Renormalisation Group approach; see for instance \cite{BGK,G2}.
We proceed as follows. Choose a subsequence $\al_{m_n}(\om)$, $n\ge 0$,
in such a way that it is strictly decreasing, and define an integer vector $\nu\in\ZZZ^d\setminus\{0\}$
to be on scale $n\ge1$ if $\al_{m_{n+1}}(\om) \le |\om\cdot\nu| < \al_{m_n}(\om)$ and on scale $n=0$ if
$|\om\cdot\nu| \ge \al_{m_0}(\om)$ (actually, for technical reasons, we shall use a smooth partition).
Roughly, the higher the scale of $\nu$ is, the smaller the small divisor associated with $\nu$.
Introducing the notation
\begin{equation} \nonumber
b=\sum_{n \ge0} b_{n},\qquad b_{\ge m} = \sum_{n \ge m} b_{n} , \qquad
b_{n} = \sum_{\nu\text{ on scale }n}{\rm e}^{\ii\nu\cdot\om t}b_{n,\nu} ,
\end{equation}
the idea is to solve the range equation ``scale by scale'', by fixing iteratively the ``components" $b_n$ in terms
of the ``corrections'' $b_{\ge n+1}$, i.e. $b_{n}=b_{n}(b_{\ge n+1})$.
Starting from scale $n=0$ we obtain from \eqref{eq:1.8}
\begin{equation}\nonumber
(\om\cdot\nu)^2 b_{0,\nu}- [ \e F(\be_0+b_0+b_{\ge 1})]_{\nu}=0,
\end{equation}
so that linearising at $b=0$ we have
\begin{equation}\nonumber
(\om\cdot\nu)^2 b_{0,\nu} - [\e \del F(\be_0)b_0]_{\nu} - [\e\del F(\be_0)b_{\ge1}]_{\nu} + O(b^2)=0 ,
\end{equation}
with $\del$ denoting derivative with respect to the argument.
Then we correct the differential operator by adding the diagonal part (with respect to the Fourier basis)
of $-\e \del F(\be_0)$; denoting $\MM_0:=\diag_\nu(\e \del F(\be_0))$ and $\NN_0:=\e\del F(\be_0)-\MM_0$,
we write
\begin{equation}\nonumber
\mathcal{D}_0(\nu) \, b_{0,\nu}=
[\NN_0b_0]_\nu + [\e \del F(\be_0)b_{\ge1}]_{\nu} + O(b^2) ,
\qquad
\mathcal{D}_0(\nu):=(\om\cdot\nu)^2\uno-\MM_{0},
\end{equation}
where $\uno$ is the $r\times r$ identity matrix. If the equation can be solved, we obtain $b_0$
expressed in terms of $b_{\ge1}$. Iterating, at the $n$-th step we write
$b=b_{\le n-1}(b_{\ge n})+b_{\ge n}$, where
\begin{equation} \nonumber
b_{\le n-1}(b_{\ge n}) := \sum_{m \le n-1} b_{m,n}(b_{\ge n}) ,
\end{equation}
with each $b_{m,n}(b_{\ge n})$ recursively expressed in terms of $b_{\ge n}$; an explicit computation
gives $b_{n-1,n}(x)=b_{n-1}(x)$ and $b_{m,n}(x) = b_{m} (b_{m+1,n}(x) + \ldots + b_{n-1,n}(x)+x)$
for $0 \le m \le n-2$. Hence, for $\nu$ on scale $n$, \eqref{eq:1.8} reduces to
\begin{equation}\nonumber
(\om\cdot\nu)^2 b_{n,\nu}- [\e F(\be_0+b_{\le n-1}(b_{\ge n})+b_{\ge n})]_{\nu}=0.
\end{equation}
We linearise at $b_{\ge n}=0$, so as to obtain
\begin{equation} \nonumber
(\om\cdot\nu)^2 b_{n,\nu} - [\e \del F(\be_0 + b_{\le n-1}(0)) (\uno + \del b_{\le n-1}(0)) \,
(b_n+b_{\ge n+1})]_{\nu} + O(b^2_{\ge n})=0 ,
\end{equation}
and add the correction $\diag_{\nn}( -\e \del F(\be_0+b_{\le n-1}(0))(\uno+\del b_{\le n-1}(0)))$
to the differential operator. In other words, the corrections to the small divisor at each step are diagonal
with respect to the Fourier indices $\nu$. Moreover they turn out to depend on $\nu$ only through
$\om\cdot\nu$, so that we can denote the correction $-\MM_{n}(\oo\cdot\nn;\e,\be_0)$.
In conclusion, at each step $n$, we have to deal with an equation of the form
\begin{equation} \label{eq:1.20}
\mathcal{D}_{n}(\nn) \, b_{n,\nn} = [ \NN_{n} b_n ]_{\nu} +
[\e \del F(\be_0 + b_{\le n-1}(0)) (\uno + \del b_{\le n-1}(0)) \, b_{\ge n+1}]_{\nu} + O(b^2_{\ge n}) ,
\end{equation}
where $\NN_{n}:=\e \del F(\be_0 + b_{\le n-1}(0)) (\uno + \del b_{\le n-1}(0)) - 
\MM_{n}(\oo\cdot\nn;\e,\be_0)$ and the operator $\mathcal{D}_{n}(\nn)$ has the form
\begin{equation} \label{eq:1.21}
\mathcal{D}_{n}(\nn) = (\om\cdot\nu)^2\uno-\MM_n(\om\cdot\nn;\e,\be_0) .
\end{equation}
By using a version of the Siegel-Bryuno bound (see Lemma \ref{lem:2.1}), we show that
the convergence of such iterative scheme would follow if one had a bound like
\begin{equation}\label{eq:1.22}
\|(\om\cdot\nu)^2\uno-\MM_n(\om\cdot\nu;\e,\be_0)\|\ge \frac{(\om\cdot\nu)^2}{2},
\end{equation}
with $\|\cdot\|$ denoting the $r$-dimensional $L^2$ operator norm.
A second order Taylor expansion of \eqref{eq:1.21} at $\om\cdot\nu=0$ gives
\begin{equation}\label{eq:1.23}
(\om\cdot\nu)^2\uno - \left( \MM_n(0;\e,\be_0)+\del \MM_n(0;\e,\be_0)(\om\cdot\nu) + 
O(\e(\om\cdot\nu)^2) \right).
\end{equation}
The time reversibility assumption in Hypothesis \ref{hyp:2} ensures that $\MM_n(\om\cdot\nn;\e,\be_0)$
is even in $\om\cdot\nn$, so that the linear term in \eqref{eq:1.23} vanishes identically and hence one
has to control $\MM_n(0;\e,\be_0)$ only. Formally (i.e. assuming the bound \eqref{eq:1.22})
one can show that, up to corrections, $\MM_n(0;\e,\be_0)=-\del_{\be_0}^2 L_n(\e,\be_0)$,
where $L_n(\e,\be_0)$ is the average of the Lagrangian computed along
the $n$-step approximate solution $b^{\le n}:= b_{\le n}(0)$.
If one could take $\be_0$ as a maximum of $L_{n}(\e,\be_0)$, the eigenvalues of $\MM_n(0;\e,\be_0)$
would be all non-positive, so  implying the bound \eqref{eq:1.22}. Obviously this is not the right way
to proceed since $L_{n}(\e,\be_0)$ is not defined for all $\be_0$ and in any case
$\be_0$ cannot be fixed at a different value at every intermediate step $n$.

The idea is to define recursively an auxiliary function $\ol{b}(\al;\e,\be_0)$, obtained
recursively from \eqref{eq:1.20} by modifing the operators \eqref{eq:1.21}
in such a way that a bound like \eqref{eq:1.22}
automatically holds. To this aim we replace recursively $\MM_n(0;\e,\be_0)$ with
$-\del_{\be_0}^2\ol{L}_n (\e,\be_0) \ol{\xi}_n(\e,\be_0)$, where $\ol{L}_n(\e,\be_0)$
is the average of the Lagrangian computed along the $(n-1)$-step approximation of the auxiliary function
and $\ol{\xi}_n$ is a suitable cut-off function.  The presence of the cut-off functions allows
the auxiliary function to be well defined for all $\be_0\in\TTT^{r}$ and thus one can pass
to the limit $n\to\io$. Of course for $\be_0$ varying in $\TTT^r$ the function $\ol{b}(\al;\e,\be_0)$
is no longer a solution of the range equation. Then we show recursively that
taking $\be_0=\ol{\be}_0(\e)$ as a maximum of $\ol{L}_{\io}(\e,\be_0)$, one has
$\MM_n(0;\e,\ol{\be}_0(\e))=-\del_{\be_0}^2\ol{L}_n (\e,\ol{\be}_0(\e))\ol{\xi}_n(\e,\ol{\be}_0(\e))$,
implying that $b(\al;\e,\ol{\be}_{0}(\e))=\ol{b}(\al;\e,\ol{\be}_0(\e))$.
Therefore both the range and the bifurcation equations are satisfied.

We note that, unlike \cite{CG1}, because of the higher dimension we cannot consider 
directly the linearisation $\MM_{n}(0;\e,\be_0)$ and then construct its second antiderivatives
with respect to $\be_0$.
On the contrary, we have to take full advantage of the variational structure of the
Hamilton equations, hence compute the Lagrangian along the solution of the range equation
and thence exploit the fact that $\MM_{n}(0;\e,\be_0)$ is formally its second derivative.
As explained above, it is not obvious that this can be achieved, because there is a priori
no control on the set of values of $\be_0$ for which all those functions are well defined:
to do this we have to introduce the auxiliary function for which everything is well defined
and eventually to prove that, by carefully choosing $\be_0$, such a function reduces
to the true solution.

%%%%%%%%%%%%%%%%%%%%%%%%%%%%%%%%%%%%%%%%%%%%%%%%%
%%%%%%%%%%%%%%%%%%%%%%%%%%%%%%%%%%%%%%%%%%%%%%%%%
\zerarcounters 
\section{Multiscale analysis}
\label{sec:2} 
%%%%%%%%%%%%%%%%%%%%%%%%%%%%%%%%%%%%%%%%%%%%%%%%%
%%%%%%%%%%%%%%%%%%%%%%%%%%%%%%%%%%%%%%%%%%%%%%%%%

We start by introducing the tree formalism; this is an efficient way to take into account the
combinatorics in the expression of the formal solution $b$ described in Section \ref{sec:1.3}.
We shall associate a numerical value with each tree in such a way that the sum
over trees of such values equals the components $b_j$ of $b$, with $j\in\{1,\ldots,r\}$
 As we shall see, many notations are very similar to those
in \cite{CG1} apart from the fact that here we are dealing with an $r$-dimensional problem
and hence we need to take into account the components $j\in\{1,\ldots,r\}$. In particular
some of the results in \cite{CG1}  do not depend on $r$ and hence hold word by word
also in the present case. Some other bounds discussed in the first
part of the next section differ from \cite{CG1} only because of the presence of $r$-dependent
constants, but they can be proved by reasoning essentially in the same way.

%%%%%%%%%%%%%%%%%%%%%%%%%%%%%%%%%%%%%%%%%%%%%%%%%
\subsection{Oriented trees}
%%%%%%%%%%%%%%%%%%%%%%%%%%%%%%%%%%%%%%%%%%%%%%%%%
 
An \emph{oriented tree} $\theta$ is a graph (that is a set of points and lines connecting them)
with no cycle,  such that all the lines are oriented toward a single point 
(\emph{root}) which has only one incident line $\ell_{\theta}$ (\emph{root line}). 
All the points in a tree except the root are called \emph{nodes}. 
The orientation of the lines in a tree induces a partial ordering 
relation ($\preceq$) between the nodes and the lines: we can 
imagine that each line carries an arrow pointing toward the root. 
Given two nodes $v$ and $w$, we write $w \prec v$ every time $v$ is along the path 
(of lines) which connects $w$ to the root. Given a node $v$,
we denote by $\p(v)$ the unique point immediately following $v$.
 
We denote by $N(\theta)$ and $L(\theta)$ the sets of nodes and 
lines in $\theta$ respectively. 
Since a line $\ell\in L(\theta)$ is uniquely identified 
by the node $v$ which it leaves, we may write $\ell = \ell_{v}$. 
We write $\ell_{w} \prec \ell_{v}$ if $w\prec v$, and $w\prec\ell=\ell_{v}$ if $w\preceq v$; 
if $\ell$ and $\ell'$ are two comparable lines, i.e. 
$\ell' \prec \ell$, we denote by $\calP(\ell,\ell')$ the 
(unique) path of lines connecting $\ell'$ to $\ell$, with $\ell$ and 
$\ell'$ not included (in particular $\calP(\ell,\ell')=\emptyset$ if $\ell'$ enters the node $\ell$ exits). 
 
Given a tree $\theta$ we call \emph{order} of $\theta$ the number $k(\theta)=|N(\theta)|=|L(\theta)|$
(for any finite set $S$ we denote by $|S|$ its cardinality).
A subset $T\subset\theta$ is a \emph{subgraph} of $\theta$
if it is formed by a set of nodes $N(T)\subseteq N(\theta)$ and
a set of lines $L(T)\subseteq L(\theta)$
in such a way that $N(T)\cup L(T)$ is connected.
If $T$ is a subgraph of $\theta$ we call \emph{order}
of $T$ the number $k(T)=|N(T)|$. We say that a line enters $T$ if it 
connects a node $v\notin N(T)$ to a node $w\in N(T)$, 
and we say  that a line exits $T$ if it connects a node $v\in N(T)$ 
to a node $w\notin N(T)$ or to the root (which is not included in $T$
in this case). Of course, if a line $\ell$ enters or exits $T$, then $\ell\notin L(T)$.

%%%%%%%%%%%%%%%%%%%%%%%%%%%%%%%%%%%%%%%%%%%%%%%%%
\subsection{Labels}
%%%%%%%%%%%%%%%%%%%%%%%%%%%%%%%%%%%%%%%%%%%%%%%%%
 
In \eqref{eq:1.4} we can write
\begin{equation} \nonumber
f(\al,\be) = \sum_{\nn\in\ZZZ^{d}} {\rm e}^{\ii\nn\cdot\al} f_{\nn}(\be) ,
\end{equation}
where $|\partial_{\be}^{s} f_{\nn}(\be)| \le s!\Phi_{0}\Phi_{1}^{s}{\rm e}^{-\xi|\nn|}$
for suitable positive constants $\Phi_{0}$ and $\Phi_{1}$,
by the analyticity assumption on $f$; note that
$f_{-\nn}(\be)=f_{\nn}(\be)$ by Hypothesis \ref{hyp:2}.
 
With each node $v\in N(\theta)$ we associate a \emph{mode} label 
$\nu_{v}\in \ZZZ^{d}$ and we denote by $s_{v}$ the number of lines entering $v$.
With each line $\ell$ we associate a \emph{momentum} $\nn_{\ell}\in \ZZZ^{d}_{*}$, 
except for the root line which can have either zero momentum or not
(i.e. $\nn_{\ell_{\theta}}\in\ZZZ^{d}$), and a pair of \emph{component labels}
$(e_{\ell},u_{\ell})\in\{1,\ldots,r\}$. We call \emph{total momentum} of 
$\theta$ the momentum associated with $\ell_{\theta}$. 
Finally, we associate with each line $\ell$ also a \emph{scale label} 
such that $n_{\ell}=-1$ if $\nn_{\ell}=\vzero$, while 
$n_{\ell}\in\ZZZ_{+}$ if $\nn_{\ell}\neq\vzero$;
note that one can have $n_{\ell}=-1$ only if $\ell$ is the root line  of $\theta$. 
We impose the following \emph{conservation law} 
\begin{equation}\label{eq:2.1} 
\nn_{\ell}=\sum_{\substack{ w \in N(\theta) \\ w\prec \ell}}\nn_{w}. 
\end{equation} 

A \emph{labelled oriented tree} is an oriented tree with labels associated
with its nodes and lines. In the following, for simplicity's sake, we shall call trees tout court
the labelled oriented trees and we shall term \emph{unlabelled tree} 
the oriented trees without labels. 
We shall say that two trees are \emph{equivalent} if they can be 
transformed into each other by continuously deforming the lines in 
such a way that these do not cross each other and also labels match. 
This provides an equivalence relation on the set of the trees. From
now on we shall call trees such equivalence classes. 

%%%%%%%%%%%%%%%%%%%%%%%%%%%%%%%%%%%%%%%%%%%%%%%%%
\subsection{Clusters and self-energy clusters}
%%%%%%%%%%%%%%%%%%%%%%%%%%%%%%%%%%%%%%%%%%%%%%%%%
 
 A \emph{cluster} $T$ on scale $n$ is a maximal subgraph 
of a tree $\theta$ such that all the lines have scales 
$n'\le n$ and there is at least a line with scale $n$. 
The lines entering the cluster $T$ and the line coming 
out from it (unique if existing at all) are called the 
\emph{external} lines of $T$. 
 
A \emph{self-energy cluster} is a cluster $T$ such that 
(i) $T$ has only one entering line $\ell'_{T}$ and 
one exiting line $\ell_{T}$, (ii) one has $\nn_{\ell_{T}}= 
\nn_{\ell'_{T}}$ and hence $\sum_{v\in N(T)}\nn_{v}=\vzero$ by \eqref{eq:2.1}.
 
For any self-energy cluster $T$, set $\calP_{T}= 
\calP(\ell_{T},\ell'_{T})$. We shall say that a self-energy cluster is on 
scale $-1$, if $N(T)=\{v\}$ with of course $\nn_{v}=\vzero$ 
(so that $\calP_{T}=\emptyset$). 
Given a self-energy cluster $T$, for all $\ell\in\calP_{T}$ one can write
\begin{equation} \label{eq:2.2}
\nu_{\ell} = \nu_{\ell}^{0}+\nu_{\ell'_{T}} , \qquad 
\nu_{\ell}^{0} := \sum_{\substack{w\in N(T) \\ w\prec \ell}} \nu_{w} .
\end{equation}
Note that the momenta of the lines along $\calP_{T}$ are the only labels
of $T$ depending on the labels outside $T$.
We say that two self-energy clusters $T_{1}$ and $T_{2}$ have the same 
\emph{structure} if setting $\nn_{\ell'_{T_{1}}}=\nn_{\ell'_{T_{2}}}=\vzero$ one has 
$T_{1}=T_{2}$. Of course this provides an equivalence relation on the 
set of all self-energy clusters; from now on we shall call 
self-energy clusters tout court such equivalence classes. 
 
%%%%%%%%%%%%%%%%%%%%%%%%%%%%%%%%%%%%%%%%%%%%%%%%%
\subsection{Renormalised trees: node factors and propagators}
%%%%%%%%%%%%%%%%%%%%%%%%%%%%%%%%%%%%%%%%%%%%%%%%%

Once we associate a numerical value with each labelled tree, it will be clear that a self-energy cluster
represents a contribution to the diagonal of the linearised vector field described in Section \ref{sec:1.3}.
Our resummation procedure implies that the self-energy clusters must appear
in the small divisor and not in the vector field. Therefore we shall consider
\emph{renormalised trees}, i.e. trees in which no self-energy 
clusters appear; analogously a \emph{renormalised subgraph} is a 
subgraph of a tree $\theta$ which does not contains any self-energy cluster. 
Denote by $\Theta_{k,\nn,j}$ the set of renormalised trees 
with order $k$, total momentum $\nn$ and $e_{\ell_{\theta}}=j$, and by 
$\gotR_{n,u,e}$ the set of renormalised self-energy clusters on scale $n$
such that $u_{\ell_{T}}=u$ and $e_{\ell_{T}'}=e$.
 
For any $\theta\in\Theta_{k,\nn,j}$ and any subgraph $S \subseteq \theta$ we associate 
with each node $v\in N(S)$ a \emph{node factor} 
\begin{equation} \label{eq:2.3}
\calF_{v}(\be_{0}) := \frac{1}{s_{v}!} \partial_{\be_{u_{\ell_{v}}}}
\Bigl( \prod_{\substack{w \in N(S) \\ \p(w)=v}} \partial_{\be_{e_{\ell_{w}}}} \Bigr)
f_{\nn_{v}}(\be_{0}). 
\end{equation} 

Introduce a partition of unity as follows.
Given a decreasing sequence $\rho_{n}$, $n\in\NNN\cup\{0\}$,
of positive numbers with $\rho_{n+1} \le \rho_{n}/2$, 
let $\chi:\RRR\to\RRR$ be a $C^{\io}$ function, non-increasing for
$x\ge0$ and non-decreasing for $x<0$, such that 
\begin{equation} \nonumber
\chi(x)=\left\{ 
\begin{aligned} 
&1,\qquad |x| \le 1/2, \\ 
&0,\qquad |x| \ge 1 , 
\end{aligned}\right.
\end{equation} 
and set $\chi_{n}(x)=\chi(x/\rho_{n})$ for $n\ge 0$ and $\chi_{-1}(x)=1$.
Set also $\Psi_{n}(x)=\chi_{n-1}(x)-\chi_{n}(x)$ for $n\ge -1$; see Figure 1 in \cite{CG1}.

Next, we introduce the sequences $\{m_{n},p_{n}\}_{n \ge 0}$, with
$m_{0}=0$ and, for all $n\ge 0$, $m_{n+1}=m_{n}+p_{n}+1$,  where
$p_{n}:=\max\{q\in\ZZZ_{+}\,:\,\al_{m_{n}}(\oo)<2\al_{m_{n}+q}(\oo)\}$,
with $\al_{m}(\oo)$ defined in (\ref{eq:1.6}).
The subsequence $\{\al_{m_{n}}(\oo)\}_{n\ge 0}$
of $\{\al_{m}(\oo)\}_{m\ge0}$ is decreasing. 
A convenient partition of unity is then obtained by choosing
$\rho_{n}=\al_{m_{n}}(\oo)/8$, which is the same choice as in \cite{CG1}.

For $n\ge0$, define formally
\begin{subequations} \label{eq:2.4} 
\begin{align}
\calG^{[n]}_{e,u}(x;\e,\be_{0}) & := \Psi_{n}(x)\left[ \left( x^{2} \uno -\MM^{[n-1]}(x;\e,\be_{0}) \right)^{-1}\right]_{e,u},
\label{eq:2.4a}\\
\MM^{[n-1]}(x;\e,\be_{0}) & := \sum_{q=-1}^{n-1}\chi_{q}(x)M^{[q]}(x;\e,\be_{0}),
\label{eq:2.4b}\\
M^{[q]}_{u,e}(x;\e,\be_{0}) & :=
\sum_{T\in\gotR_{q,u,e}}\e^{k(T)}\Val_{T}(x;\e,\be_{0}),
\label{eq:2.4c} \\
\Val_{T}(x;\e,\be_{0}) & := \Biggl(\prod_{v\in N(T)}\calF_{v}(\be_{0}) \Biggr)
\Biggl(\prod_{\ell\in L(T)}\calG^{[n_{\ell}]}_{e_{\ell},u_{\ell}} (\oo\cdot\nn_{\ell};\e,\be_{0})\Biggr),
\label{eq:2.4d}
\end{align}
\end{subequations} 
where $x=\oo\cdot\nn_{\ell_{T}'}$ and
$\Val_{T}(x;\e,\be_{0})$ is called the \emph{renormalised value} of $T$.
Here and henceforth, the sums and the products over empty sets have to be
considered as 0 and $1$, respectively.
Set $\MM=\{\MM^{[n]}(x;\e,\be_{0})\}_{n\ge-1}$. 
We call \emph{self-energies} the $r\times r$ matrices $\MM^{[n]}(x;\e,\be_{0})$. 

Then we associate with each line $\ell\in L(\theta)$ a \emph{propagator} $\calG_{\ell}$ by setting
\begin{equation} \nonumber
\calG_{\ell} = \begin{cases}
\calG^{[n_{\ell}]}_{e_{\ell},u_{\ell}}(\oo\cdot\nn_{\ell};\e,\be_0) , & \qquad n_{\ell} \ge 0 , \\
\de_{e_{\ell},u_{\ell}} , & \qquad n_{\ell}=-1 .
\end{cases}
\end{equation}
where $\de_{e,u}$ is the Kronecker symbol.
Recall that $n_{\ell}=-1$ is possible only if $\ell$ is the root line.

Note that, in defining the renormalised value in \eqref{eq:2.4a}
as $\Val_{T}(\oo\cdot\nn_{\ell_{T}'};\e,\be_0)$,
we have used that, by construction, the propagators of the lines $\ell\in\calP_{T}$
depend on $\oo\cdot\nn_{\ell_{T}'}$, while for the propagators of the lines
$\ell\in L(T)\setminus\calP_{T}$ one has $\nn_{\ell}=\nn_{\ell}^{0}$,
with the notation \eqref{eq:2.2}.

%%%%%%%%%%%%%%%%%%%%%%%%%%%%%%%%%%%%%%%%%%%%%%%%%
\subsection{Resummed series}
%%%%%%%%%%%%%%%%%%%%%%%%%%%%%%%%%%%%%%%%%%%%%%%%%

For any subgraph $S$ of any  
$\theta\in \Theta_{k,\nn,j}$ define the \emph{renormalised value} of $S$ as 
\begin{equation} \nonumber
\Val(S;\e,\be_{0}) := \Biggl(\prod_{v\in N(S)}\calF_{v}(\be_{0}) 
\Biggr) \Biggl(\prod_{\ell\in L(S)}\calG_{\ell} \Biggr). 
\end{equation} 

Set 
\begin{subequations} \label{eq:2.5} 
\begin{align}
b_{\nn j}^{[k]}(\e,\be_{0}) & := \sum_{\theta\in \Theta_{k,\nn,j}}\Val (\theta;\e,\be_{0}), 
\quad \nn \neq \vzero ,
\label{eq:2.5a} \\
G^{[k]}_{j}(\e,\be_{0}) & := \sum_{\theta\in \Theta_{k+1,\vzero,j}}\Val(\theta; \e,\be_{0}), 
\label{eq:2.5b}
\end{align}
\end{subequations} 
and define formally the functions $b(\al;\e,\be_{0})$ and $G(\e,\be_{0})$ with components
\begin{subequations} \label{eq:2.6} 
\begin{align}
b_{j}(\al;\e,\be_{0}) & := \sum_{k\ge1}\e^{k}\sum_{\nn\in\ZZZ_{*}^{d}} 
{\rm e}^{\ii\nn\cdot\al} b_{\nn j}^{[k]}(\e,\be_{0}), 
\label{eq:2.6a} \\
G_{j}(\e,\be_{0}) &  := \sum_{k\ge 0}\e^{k+1} G^{[k]}_{j}(\e,\be_{0}). 
\label{eq:2.6b}
\end{align}
\end{subequations} 
We call (\ref{eq:2.6}) the \emph{resummed series}. Note that formally
$G(\e,\be_0)=[\e \del_{\be}f(\cdot,\be_0+b(\cdot;\e,\be_0))]_{0}$ and $b=b(\al;\e,\be_0)$ is
the formal solution referred to in Section \ref{sec:1.3}.
 
%%%%%%%%%%%%%%%%%%%%%%%%%%%%%%%%%%%%%%%%%%%%%%%%%
\subsection{Siegel-Bryuno bounds}
%%%%%%%%%%%%%%%%%%%%%%%%%%%%%%%%%%%%%%%%%%%%%%%%%

For $\theta\in\Theta_{k,\nn,j}$, let $\gotN_{n}(\theta)$ be the number 
of lines on scale $\ge n$ in $\theta$, and set 
\begin{equation} \nonumber 
K(\theta):=\sum_{v\in N(\theta)}|\nn_{v}|. 
\end{equation} 
More generally, for any renormalised subgraph $T$ of any tree $\theta$ 
call $\gotN_{n}(T)$ the number of lines on scale $\ge n$ in $T$, and set 
\begin{equation} \nonumber 
K(T):=\sum_{v\in N(T)}|\nn_{v}|. 
\end{equation} 

The forthcoming two results are identical to Lemmas 4.1 and 4.2 in \cite{CG1},
respectively, since they do not depend on $r$.

%%%%%%%%%%%%%%%%%%%%%%%%%%%%%%%%%%%%%%%%%%%%%%%%%
\begin{lemma}\label{lem:2.1} 
For any $\theta\!\in\!\Theta_{k,\nn,j}$ 
such that $\Val(\theta;\e,\be_{0}) \!\neq \!0$ one has $\gotN_{n}(\theta) 
\!\le\! 2^{-(m_{n}-2)}K(\theta)$ for all $n\ge0$. 
\end{lemma} 
%%%%%%%%%%%%%%%%%%%%%%%%%%%%%%%%%%%%%%%%%%%%%%%%%
 
The proof is as in \cite{CG1}, Appendix A.
 
%%%%%%%%%%%%%%%%%%%%%%%%%%%%%%%%%%%%%%%%%%%%%%%%%
\begin{lemma}\label{lem:2.2} 
For any $T\!\in\!\gotR_{n,u,e}$ such that $\Val_{T}(x;\e,\be_{0})\! \neq \!0$, 
one has $K(T)\!\ge \!2^{m_{n}-1}$ and $\gotN_{p}(T)\!\le
\!2^{-(m_{p}-2)} K(T)$ for all $0\le p\le n$. 
\end{lemma} 
%%%%%%%%%%%%%%%%%%%%%%%%%%%%%%%%%%%%%%%%%%%%%%%%%

The proof is as in \cite{CG1}, Appendix B.
  
%%%%%%%%%%%%%%%%%%%%%%%%%%%%%%%%%%%%%%%%%%%%%%%%%
%%%%%%%%%%%%%%%%%%%%%%%%%%%%%%%%%%%%%%%%%%%%%%%%%
\zerarcounters 
\section{Convergence of the resummed series}
\label{sec:3} 
%%%%%%%%%%%%%%%%%%%%%%%%%%%%%%%%%%%%%%%%%%%%%%%%%
%%%%%%%%%%%%%%%%%%%%%%%%%%%%%%%%%%%%%%%%%%%%%%%%%

We first show that, under the assumption that the propagators satisfy suitable bounds
(that we call property 1), the convergence of the resummed series would follow
by standard arguments of multiscale analysis. As we shall see, two fundamental
ingredients of the proof will be a cancellation mechanism, which strongly relies
on the symmetry in Hypothesis \ref{hyp:2} (Lemma \ref{lem:3.5}), and some remarkable
identities which relate the self-energies to the averaged Lagrangian (Lemma \ref{lem:3.8}).

Next, to remove the undesired assumption, we proceed as follows. We modify the propagators
by replacing the self-energies $\MM$ with new matrices $\ol{\MM}^{\xi}$, which
still satisfy the symmetry of Lemma \ref{lem:3.5}: by construction, this will imply
automatically that $\ol{\MM}^{\xi}$ satisfies property 1. Then we check that on a suitable curve
$\be_{0}=\ol{\be}_{0}(\e)$ one has $\ol{\MM}^{\xi}=\MM$ and hence the range
equation is satisfied. Finally, by further exploiting the identites of Lemma \ref{lem:3.8},
we show that, on such a curve, also the bifurcation equation is solved.

%%%%%%%%%%%%%%%%%%%%%%%%%%%%%%%%%%%%%%%%%%%%%%%%%
\subsection{Formal analysis} \label{sec:3.1}
%%%%%%%%%%%%%%%%%%%%%%%%%%%%%%%%%%%%%%%%%%%%%%%%%

We first assume that the propagators $\calG_{\ell}$ are bounded proportionally to 
$(\oo\cdot\nn_{\ell})^{-2}$ and show that, under such an assumption, the convergence of the series
is easily checked. Then, we shall prove that the assumption makes sense for a suitable choice of $\be_0$.
Let us denote by $\|\cdot\|$ the $r$-dimensional $L^2$ operator norm.

%%%%%%%%%%%%%%%%%%%%%%%%%%%%%%%%%%%%%%%%%%%%%%%%%
\begin{defi}\label{defi:3.1} 
We shall say that $\MM$ satisfies \emph{property 1-$p$} if for $-1\le n <p$ one has 
\begin{equation} \nonumber
\Psi_{n+1}(x) \left\|x^{2}\uno-\MM^{[n]}(x;\e,\be_{0})\right\|\ge \Psi_{n+1}(x) \, x^{2}/2 .
\end{equation} 
\end{defi} 
%%%%%%%%%%%%%%%%%%%%%%%%%%%%%%%%%%%%%%%%%%%%%%%%%

Denote by $\Theta^{\le p}_{k,\nn,j}$ the set of renormalised trees whose lines are
on scale $\le p$ and set
\begin{equation} \label{eq:3.1}
b^{\le p}_{j}(\al;\e,\be_{0}) := \sum_{k\ge1}\e^{k}\sum_{\nn\in\ZZZ_{*}^{d}} 
{\rm e}^{\ii\nn\cdot\al}
\sum_{\theta\in \Theta^{\le p}_{k,\nn,j}}\Val (\theta;\e,\be_{0}) .
\end{equation}
Call $b^{\le p}(\al;\e,\be_0)$ the function with components \eqref{eq:3.1}. 

%%%%%%%%%%%%%%%%%%%%%%%%%%%%%%%%%%%%%%%%%%%%%%%%%
\begin{lemma}\label{lem:3.2} 
Assume $\MM$ to satisfy property 1-$p$.
There are two positive constants $B_0$ and $B_1$ such that
for all $(k,\nn,j) \in \NNN \times \ZZZ^{d} \times \{1,\ldots,r\}$ and for any tree
$\theta\in\Theta^{\le p}_{k,\nn,j}$ one has $|\Val(\theta;\e,\be_{0}))|
\le B_{0}B_{1}^{k}{\rm e}^{-\xi_{1}|\nn|}$, with $B_0,B_1$ and $\xi_{1}<\xi$, independent of $p$.
\end{lemma} 
%%%%%%%%%%%%%%%%%%%%%%%%%%%%%%%%%%%%%%%%%%%%%%%%%

The proof is standard; see for instance \cite{CG1}, Appendix C. Note that the constants
$B_0$ and $B_1$ increase with $r$.

Lemma \ref{lem:3.2} implies immediately the following result.

%%%%%%%%%%%%%%%%%%%%%%%%%%%%%%%%%%%%%%%%%%%%%%%%%
\begin{lemma}\label{lem:3.3} 
Assume $\MM$ to satisfy property 1-$p$.
Then the series (\ref{eq:3.1}) converge for $\e$ small enough.
Moreover $b^{\le p}(\al;\e,\be_{0})$ is analytic in $\al\in\TTT^{d}_{\x_{2}}$,
with $\xi_{2}<\xi_{1}$, uniformly in $p$.
\end{lemma}
%%%%%%%%%%%%%%%%%%%%%%%%%%%%%%%%%%%%%%%%%%%%%%%%%

%%%%%%%%%%%%%%%%%%%%%%%%%%%%%%%%%%%%%%%%%%%%%%%%%
\begin{lemma}\label{lem:3.4} 
Assume $\MM$ to satisfy property 1-$p$. Then for any 
$0\le n\le p$ the self-energies are well defined and one has 
\begin{equation} \nonumber 
\left| \partial^{j}_{x}M^{[n]}_{e,u}(x;\e,\be_{0}) \right|
\le |\e|^{2}C_{0} {\rm e}^{-C_{1} 2^{m_{n}}},  \qquad j=0,1,2, 
\end{equation} 
for suitable constants $C_{0}$ and $C_{1}$, independent of $p$.
\end{lemma} 
%%%%%%%%%%%%%%%%%%%%%%%%%%%%%%%%%%%%%%%%%%%%%%%%%
 
The proof is as in \cite{CG1}, Appendix E.

By writing
\begin{equation} \nonumber
\MM^{[n]}(x;\e,\be_{0}) = \MM^{[n]}(0;\e,\be_{0}) + x \,
\partial_{x} \MM^{[n]}(0;\e,\be_{0}) + x^2
\int_{0}^{1} {\rm d}\tau \, \left( 1-\tau \right) \partial_{x}^{2}
\MM^{[n]}(\tau x;\e,\be_{0}) ,
\end{equation}
one checks easily that, if $\MM$ satisfies property 1-$p$, $\partial_{x}^{j}M^{[n]}_{e,u}(\tau x;\e,\be_{0})$
admits the same bounds as in Lemma \ref{lem:3.4}, for $0\le n\le p$, $j=0,1,2$ and $\tau\in[0,1]$. This
implies that
\begin{equation} \nonumber
\left\| \MM^{[n]}(x;\e,\be_{0})-\MM^{[n]}(0;\e,\be_{0}) -x \, 
\partial_{x} \MM^{[n]}(0;\e,\be_{0}) \right\| \le C_{2} |\e|^{2} x^2
\end{equation}
for some positive constant $C_{2}$ independent of $p$.

%%%%%%%%%%%%%%%%%%%%%%%%%%%%%%%%%%%%%%%%%%%%%%%%%
\begin{lemma}\label{lem:3.5} 
Assume $\MM$ to satisfy property 1-$p$. Then one has
$\del_{x} \MM^{[n]}(0;\e,\be_{0})=0$ for all $-1 \le n\le p$. 
\end{lemma} 
%%%%%%%%%%%%%%%%%%%%%%%%%%%%%%%%%%%%%%%%%%%%%%%%%
 
%%%%%%%%%%%%%%%%%%%%%%%%%%%%%%%%%%%%%%%%%%%%%%%%%
\prova 
One proves that $\MM^{[n]}(x;\e,\be_{0})=(\MM^{[n]}(-x;\e,\be_{0}))^T=
(\MM^{[n]}(x;\e,\be_{0}))^{\dagger}$ by induction as in \cite{CFG}, \S 4; see also \cite{GG}.
The key remark is that, if we replace the momentum $\nn$ of the line $\ell_{T}'$
with $-\nn$ and then reverse the orientation of the lines in $\calP_T\cup\{\ell_T,\ell'_T\}$,
% exchange the entering line $\ell_{T}'$ with the exiting line $\ell_{T}$,
then the momenta of the lines $\ell\in L(T)\setminus \calP_{T}$ do not change,
whereas the momenta of the lines $\ell\in\calP_{T}$ change sign.
In particular one has $\MM^{[n]}(-x;\e,\be_{0})=(\MM^{[n]}(x;\e,\be_{0}))^*$.
Moreover Hypothesis \ref{hyp:2} yields that the node factors \eqref{eq:2.3} are real and
this implies inductively that also the propagators are real, so that
$(\MM^{[n]}(x;\e,\be_{0}))^{*}=\MM^{[n]}(x;\e,\be_{0})$.
Therefore $\MM^{[n]}(-x;\e,\be_{0}) = \MM^{[n]}(x;\e,\be_{0})$.\EP
%%%%%%%%%%%%%%%%%%%%%%%%%%%%%%%%%%%%%%%%%%%%%%%%%

%%%%%%%%%%%%%%%%%%%%%%%%%%%%%%%%%%%%%%%%%%%%%%%%%
\begin{rmk} \label{rmk:3.6}
\emph{
Lemma \ref{lem:3.5}, which is crucial in order to obtain the final result,
is the only point where the reversibility condition in Hypothesis \ref{hyp:2} is used.
For $r=1$ we do not need such an assumption, because in that case
$\MM^{[n]}(x;\e,\be_0)$ is number-valued (instead of matrix-valued) and hence the identity
$\MM^{[n]}(x;\e,\be_0)=(\MM^{[n]}(-x;\e,\be_0))^{T}$ is enough to get the analogous result.
Of course the symmetry $\MM^{[n]}(-x;\e,\be_{0}) = \MM^{[n]}(x;\e,\be_{0})$,
for which the parity assumption seems necessary, 
is much stronger than the desired result $\del_{x} \MM^{[n]}(0;\e,\be_{0})=0$.
}
\end{rmk}
%%%%%%%%%%%%%%%%%%%%%%%%%%%%%%%%%%%%%%%%%%%%%%%%%

%%%%%%%%%%%%%%%%%%%%%%%%%%%%%%%%%%%%%%%%%%%%%%%%%
\begin{lemma} \label{lem:3.7}
Assume $\MM$ to satisfy property 1-$p$.
Then, for all $\be_0 \in \TTT^{r}$, $b^{\le p}(\al;\e,\be_0)$ solves
the range equation \eqref{eq:1.7a}
up to corrections bounded by $C_{3}{\rm e}^{-\xi' 2^{m_{p+1}}}$,
for some positive constants $C_{3}$ and $\xi'$ independent of $p$.
\end{lemma}
%%%%%%%%%%%%%%%%%%%%%%%%%%%%%%%%%%%%%%%%%%%%%%%%%

%%%%%%%%%%%%%%%%%%%%%%%%%%%%%%%%%%%%%%%%%%%%%%%%%
\prova
We want to show that, by shortening $b^{\le p}=b^{\le p}(\al;\e,\be_0)$,
\begin{equation} \nonumber
(\oo\cdot\nn)^{-2} \left[ \e \del_\be f(\al,\be_{0}+b^{\le p} ) \right]_{\nn} = 
[b^{\le p}]_{\nn} + O({\rm e}^{-\xi' 2^{m_{p+1}}})
\end{equation}
for some $\xi'>0$. For any $n\le p$,
denote by $\Theta^{\le p}_{k,\nn,j}(n)$ the set of trees $\theta\in \Theta^{\le p}_{k,\nn,j}$
such that the root line is on scale $n$ and set
\begin{equation} \nonumber
b^{n}_{\nu} := \sum_{k\ge 1} \e^{k}
\sum_{\theta\in \Theta^{\le p}_{k,\nn,j}(n)}\Val (\theta;\e,\be_{0}) .
\end{equation}
We can write
\begin{equation} \label{eq:3.2} 
\begin{aligned} 
& [\e \del_\be f(\al,\be_{0}+b^{\le p})]_{\nn}
= \sum_{n=0}^{\io} \Psi_{n}(\oo\cdot\nn) [\e \del_\be f(\al,\be_{0}+b^{\le p})]_{\nn} \\ 
& \quad = \sum_{n=0}^{p} \Psi_{n}(\oo\cdot\nn) [\e \del_\be f(\al,\be_{0}+b^{\le p})]_{\nn} +
\sum_{n=p+1}^{\io} \Psi_{n}(\oo\cdot\nn) [\e \del_\be f(\al,\be_{0}+b^{\le p})]_{\nn} ,
\end{aligned}
\end{equation}
with the last sum not vanishing only if $\Psi_{n}(\oo\cdot\nn) \neq0$ for
some $n\ge p+1$. In turn this requires $|\nn| \ge 2^{m_{p+1}}$, so that, by relying
on Lemma \ref{lem:3.2}, for such $\nn$ we can bound
\begin{equation} \label{eq:3.3}
c_{\nn}:=(\oo\cdot\nn)^{-2} \Bigl| \sum_{n=p+1}^{\io} \Psi_{n}(\oo\cdot\nn)
[\e \del_\be f(\al,\be_{0}+b^{\le p})]_{\nn} \Bigr|
\le C \, {\rm e}^{-\xi' |\nn|} {\rm e}^{-\xi' 2^{m_{p+1}}} .\
\end{equation}
for some positive constants $C$ and $\xi'$.
The first sum in the second line of \eqref{eq:3.2} becomes
\begin{equation} \label{eq:3.4} 
\sum_{n=0}^{p} \left((\oo\cdot\nn)^{2}-\MM^{[n-1]}(\oo\cdot\nn; \e,\be_{0})\right)
\sum_{k\ge1} \e^{k}\sum_{\theta\in\ol{\Theta}^{\le p}_{k,\nn,j}(n)}  \Val(\theta;\e,\be_{0}) ,
\end{equation} 
where $\ol{\Theta}^{\le p}_{k,\nn,j}(n)$ differs from $\Theta^{\le p}_{k,\nn,j}(n)$ 
as it contains also trees $\theta$ which have one self-energy cluster with exiting line $\ell_{\theta}$. 
If we separate the trees containing such self-energy cluster from the others, \eqref{eq:3.4} gives
\begin{equation} \nonumber 
\begin{aligned} 
& \sum_{n=0}^{p}\left((\oo\cdot\nn)^{2}-\MM^{[m-1]}(\oo\cdot\nn; \e,\be_{0})\right)
b_{\nn}^{n} + \sum_{m=0}^{p} \sum_{q=-1}^{m-1}  M^{[q]}(\oo\cdot\nn;\e,\be_{0})
\sum_{n= q+1}^{p}\Psi_{n}(\oo\cdot\nn) \, b_{\nn}^{m} \\ 
& \qquad = \sum_{n=0}^{p} \left((\oo\cdot\nn)^{2}-\MM^{[n-1]}(\oo\cdot\nn; \e,\be_{0})\right) b_{\nn}^{n} 
+ \sum_{n=0}^{p} \MM^{[n-1]}(\oo\cdot\nn;\e,\be_{0})  \, b_{\nn}^{n} , 
\end{aligned} 
\end{equation} 
so that 
\begin{equation} \nonumber
(\oo\cdot\nn)^{-2} \left[ \e \del_\be f(\al,\be_{0}+b^{\le p} ) \right]_{\nn} =
c_{\nn} + \sum_{n=0}^{p} b_{\nn}^{n} = c_{\nn} + [b^{\le p}]_{\nn} 
\end{equation}
where $c_{\nn}$ vanishes for $|\nn|<2^{m_{p+1}}$ and satisfies \eqref{eq:3.3} for $|\nn|\ge 2^{m_{p+1}}$.
\EP
%%%%%%%%%%%%%%%%%%%%%%%%%%%%%%%%%%%%%%%%%%%%%%%%%

Let us denote
\begin{equation}\label{eq:3.5}
G^{\le p}(\e,\be_0):= [ \e \del_\be f(\al,\be_{0}+b^{\le p}(\al;\e,\be_0))]_{0}
\end{equation}
and
\begin{equation} \label{eq:3.6}
L^{\le p}(\e,\be_0):= \Bigl[ -\frac{1}{2}(\om\cdot\del_\al b^{\le p}(\al;\e,\be_0))^2+
\e f(\al,\be_0+b^{\le p}(\al;\e,\be_0))\Bigr]_{0} ,
\end{equation}
which are the averages of the vector field and of the Lagrangian respectively,
both computed along the approximate solution $\be_0+b^{\le p}(\al;\e,\be_0)$.

The variational structure of the Hamilton equations implies that, if we could solve
the range equation, then the left hand side of the bifurcation equation \eqref{eq:1.7b}
would be the derivative with respect of $\be_0$ of the average of the Lagrangian
computed along the solution. A variational identity of this kind holds,
up to some correction, for the approximate solution; indeed we have the following result.

%%%%%%%%%%%%%%%%%%%%%%%%%%%%%%%%%%%%%%%%%%%%%%%%%
\begin{lemma} \label{lem:3.8}
Assume $\MM$ to satisfy property 1-$p$. 
Then \eqref{eq:3.5} and \eqref{eq:3.6} are
well defined and $C^\io$ in both $\e$ and $\be_0$.
Moreover one has for all $-1 \le n \le p$
\begin{equation} \nonumber
\MM^{[n]}(0;\e,\be_0) = \del_{\be_0} G^{\le n}(\e,\be_0) + c_n(\e,\be_0)
, \qquad
G^{\le n}(\e,\be_0) = \del_{\be_0} L^{\le n}(\e,\be_0) + c'_n(\e,\be_0) ,
\end{equation}
for suitable $c_n(\e,\be_0),c'_n(\e,\be_0)=O({\rm e}^{-\x' 2^{m_{n+1}}})$,
with $\x'$ as in Lemma \ref{lem:3.7}.
\end{lemma}
%%%%%%%%%%%%%%%%%%%%%%%%%%%%%%%%%%%%%%%%%%%%%%%%%

%%%%%%%%%%%%%%%%%%%%%%%%%%%%%%%%%%%%%%%%%%%%%%%%%
\prova
The first identity can be proved as Lemma 4.8 in \cite{CG2}.
The second one is easier: integrating by parts and using Lemma \ref{lem:3.7} and
the fact that $[b^{\le p}(\cdot;\e,\be_0)]_0=0$, we obtain,
shortening again $b^{\le p}=b^{\le p}(\al;\e,\be_0)$,
\begin{equation} \nonumber
\begin{aligned}
\del_{\be_0} L^{\le n} & = \frac{1}{(2\p)^{d}} \int_{\TTT^d}\!\!\! {\rm d}\al \Big(
\e \del_{\be} f(\al,\be_0+b^{\le n})(1+\del_{\be_0} b^{\le n}) - (\om\cdot\del_\al b^{\le n})
(\om\cdot\del_\al \del_{\be_0}b^{\le n})\Big) \\
& = \frac{1}{(2\p)^{d}} \int_{\TTT^d}{\rm d}\al \Big(
\e \del_{\be} f(\al,\be_0+b^{\le n})(1 + \del_{\be_0}b^{\le n})
+ ( (\om\cdot\del_\al)^2 b^{\le n}) \, \del_{\be_0}b^{\le n} \Big) \\
& = \frac{1}{(2\p)^{d}} \int_{\TTT^d}{\rm d}\al \Big( \e f(\al,\be_0+b^{\le n}) \Big)
+ c''_n(\e,\be_0) ,
\end{aligned}
\end{equation}
with a suitable $c_{n}''(\e,\be_0)=O({\rm e}^{-\x' 2^{m_{n+1}}})$.
\EP
%%%%%%%%%%%%%%%%%%%%%%%%%%%%%%%%%%%%%%%%%%%%%%%%%

%%%%%%%%%%%%%%%%%%%%%%%%%%%%%%%%%%%%%%%%%%%%%%%%%
\begin{defi} \label{defi:3.9} 
We shall say that $\MM$ satisfies \emph{property 1} if for all $n\ge -1$ one has 
\begin{equation} \nonumber 
\Psi_{n+1}(x)\left\|x^{2}\uno-\MM^{[n]}(x;\e,\be_{0})\right\|\ge \Psi_{n+1}(x)\,x^{2}/2.
\end{equation} 
\end{defi} 
%%%%%%%%%%%%%%%%%%%%%%%%%%%%%%%%%%%%%%%%%%%%%%%%%

%%%%%%%%%%%%%%%%%%%%%%%%%%%%%%%%%%%%%%%%%%%%%%%%%
\begin{rmk} \label{rmk:3.10}
\emph{
Assuming property 1 means assuming the bounds of property 1-$p$ for all $p\ge 0$.
The reason why we have introduced first property 1-$p$, and not directly property 1, is that
property 1-$p$ is what is needed to prove by induction that property 1 holds.
Indeed, in Section \ref{sec:3.2} we shall introduce new matrices $\ol{\MM}^{\xi}$ and we shall prove that
they satisfy property 1. This will be achieved by induction, by assuming that they satisfy
property 1-$p$ and hence showing that property 1-$(p+1)$ follows as well: to this aim
we shall need the lemmas above which require property 1-$p$ to be satisfied.
}
\end{rmk}
%%%%%%%%%%%%%%%%%%%%%%%%%%%%%%%%%%%%%%%%%%%%%%%%%

If property 1 holds, the approximate solution $b^{\le p}(\al;\e,\be_0)$
is a ``scale $p$ truncation" of the solution to the range equation. Indeed we have the following result.

%%%%%%%%%%%%%%%%%%%%%%%%%%%%%%%%%%%%%%%%%%%%%%%%%
\begin{lemma}\label{lem:3.11} 
Assume $\MM$ to satisfy property 1. Then for $\e$ small enough 
the function $b(\al;\e,\be_{0})$ in (\ref{eq:2.6a}), with the coefficients given in (\ref{eq:2.5}), 
is analytic in $\al\in\TTT^{d}_{\x_{2}}$, with $\xi_{2}$ as in Lemma \ref{lem:3.3},
and solves the equation (\ref{eq:1.7a}). 
\end{lemma} 
%%%%%%%%%%%%%%%%%%%%%%%%%%%%%%%%%%%%%%%%%%%%%%%%%

The proof follows from Lemmas \ref{lem:3.3} and \ref{lem:3.7}, by taking the limit $p\to\io$.

Define formally 
\begin{equation}\label{eq:3.7} 
\MM^{[\io]}(x;\e,\be_{0}) := \lim_{n\to\io}\MM^{[n]}(x;\e,\be_{0}), \qquad  
L^{[\io]}(\e,\be_0) := \lim_{n\to\io}L^{\le n}(\e,\be_0) .
\end{equation} 
Lemma \ref{lem:3.4} yields that, under the assumption that $\MM$ satisfy property 1,
the functions \eqref{eq:3.7} are well defined and smooth. Indeed the following result holds.

%%%%%%%%%%%%%%%%%%%%%%%%%%%%%%%%%%%%%%%%%%%%%%%%%
\begin{lemma}\label{lem:3.12} 
Assume $\MM$ to satisfy property 1. 
Then the function $G(\e,\be_{0})$ and the self-energies $\MM$, as well as the
limits \eqref{eq:3.7}, are $C^{\io}$ in both $\e$ and $\be_{0}$. 
\end{lemma} 
%%%%%%%%%%%%%%%%%%%%%%%%%%%%%%%%%%%%%%%%%%%%%%%%%
 
%%%%%%%%%%%%%%%%%%%%%%%%%%%%%%%%%%%%%%%%%%%%%%%%%
\begin{lemma}\label{lem:3.13} 
Assume $\MM$ to satisfy property 1. Then the implicit function 
equation $G(\e,\be_{0})=0$ admits a solution $\be_{0}=\be_{0}(\e)$ such that one has
$\partial_{\be_{0}}G(\e,\be_{0}(\e))\le 0$. 
\end{lemma} 
%%%%%%%%%%%%%%%%%%%%%%%%%%%%%%%%%%%%%%%%%%%%%%%%%
 
%%%%%%%%%%%%%%%%%%%%%%%%%%%%%%%%%%%%%%%%%%%%%%%%%
\prova 
Property 1 ensures that we can apply Lemma \ref{lem:3.8} with $n\to\io$,
so as to obtain $\partial_{\be_{0}}G(\e,\be_0)=\del_{\be_0}^2L^{[\io]}(\e,\be_0)$.
Since $L^{[\io]}(\e,\be_0)$ is a smooth function on the torus $\TTT^{r}$,
it admits at least one maximum.
\EP
%%%%%%%%%%%%%%%%%%%%%%%%%%%%%%%%%%%%%%%%%%%%%%%%%

Lemmas \ref{lem:3.11} and \ref{lem:3.13} give the following result.
 
%%%%%%%%%%%%%%%%%%%%%%%%%%%%%%%%%%%%%%%%%%%%%%%%%
\begin{lemma}\label{lem:3.14}
Assume $\MM$ to satisfy property 1 and let $\be(\e)$ be as in Lemma \ref{lem:3.13}.
Then the function $\be(t)=\be_{0}(\e)+b(\al;\e,\be_{0}(\e))$ solves \eqref{eq:1.5}. 
\end{lemma}
%%%%%%%%%%%%%%%%%%%%%%%%%%%%%%%%%%%%%%%%%%%%%%%%%
 
%%%%%%%%%%%%%%%%%%%%%%%%%%%%%%%%%%%%%%%%%%%%%%%%%
\subsection{Rigorous analysis: phase locking} \label{sec:3.2}
%%%%%%%%%%%%%%%%%%%%%%%%%%%%%%%%%%%%%%%%%%%%%%%%%

We introduce an auxiliary function $\ol{b}(\al;\e,\be_0)$ obtained from the tree expansion 
by replacing the propagators \eqref{eq:2.4} in such a way that the new propagators
$\ol{\calG}_{\ell}$ are bounded proportionally to $(\oo\cdot\nn_{\ell})^{-2}$.
%and hence the formal results in Section \ref{sec:3.1} apply.

For all $n\ge 0$, define the $C^{\io}$ non-increasing functions $\x_{n}$ such that 
\begin{equation} \nonumber
\x_{n}(x_1,\ldots,x_r)=\left\{ 
\begin{aligned} 
& 0, \qquad \exists i \in\{1,\ldots,r\} \hbox{ s.t. } x_i\ge \frac{\al^{2}_{m_{n+1}}(\oo)}{2^{11}},\\ 
& 1, \qquad x_{i}\le \frac{\al_{m_{n+1}}^{2}(\oo)}{2^{12}} \quad \forall i=1,\ldots,r , 
\end{aligned} \right. 
\end{equation} 
and set $\x_{-1}(x_1,\ldots,x_r)=1$.

Define, for $n\ge0$, the modified propagators as
$\ol{\calG}_{\ell}=\ol{\calG}^{[n_{\ell}]}_{e_{\ell},u_{\ell}}(\oo\cdot\nn_{\ell};\e,\be_{0})$, with
\begin{equation} \nonumber
\ol{\calG}^{[n]}(x;\e,\be_{0})=\Psi_{n}(x) 
\left( x^{2}\uno-\ol{\MM}^{[n-1]}(x;\e,\be_{0}) \, \x_{n-1}(\ol{\la}^{[n-1]}(\e,\be_{0})) \right)^{-1}, 
\end{equation} 
where $\ol{\MM}^{[-1]}(x;\e,\be_{0}) = \del_{\be_{0}}^{2} \ol{L}^{[-1]}(\e,\be_0)$, with
$\ol{L}^{[-1]}(\e,\be_0)=\e f_{\vzero}(\be_{0})$,
while for $n\ge0$ one proceeds recursively as follows:
\begin{enumerate}
\item set
\begin{equation} \label{eq:3.8} 
\ol{b}^{[n]}_{j}(\al;\e,\be_{0}) := \sum_{k\ge1}\e^{k}\sum_{\nn\in\ZZZ_{*}^{d}} 
{\rm e}^{\ii\nn\cdot\al}
\sum_{\theta\in \Theta^{\le n}_{k,\nn,j}}\ol{\Val}(\theta;\e,\be_{0}) ,
\end{equation}
where
\begin{equation} \nonumber
\ol{\Val}(\theta;\e,\be_{0}) := 
\Biggl(\prod_{v\in N(\theta)}\calF_{v}(\be_{0})\Biggr)\Biggl(\prod_{\ell\in L(\theta)} 
\ol{\calG}_{\ell} \Biggr) ,
\end{equation} 
and call $\ol{b}^{[n]}(\al;\e,\be_0)$ the function with components \eqref{eq:3.8};
\item define
\begin{equation} \nonumber
\ol{L}^{[n]}(\e,\be_0):= \Bigl[ -\frac{1}{2}(\om\cdot\del_\al \ol{b}^{[n]}(\al;\e,\be_0))^2+
\e f(\al,\be_0+\ol{b}^{[n]}(\al;\e,\be_0))\Bigr]_{0} ;
\end{equation}
\item consider any matrices $\ol{c}_{n}(\e,\be_0)$ and $\ol{\gotR}^{[n]}(x;\e,\be_0)$,
satisfying the constraints
\begin{subequations} \label{eq:3.9} 
\begin{align}
& \left\| \ol{c}_n(\e,\be_0) \right\| \le K |\e|^{2} {\rm e}^{-\x' 2^{m_{n+1}}} ,
\label{eq:3.9a} \\
& \left\| \ol{\gotR}^{[n]}(x;\e,\be_0) - \ol{\gotR}^{[n-1]}(x;\e,\be_0) \right\| \le
K |\e|^{2} x^{2} {\rm e}^{-\x' 2^{m_{n}}} ,
\label{eq:3.9b}
\end{align}
\end{subequations}
for some positive constant $K$, with $\ol{\gotR}^{[-1]}(x;\e,\be_0)=0$ and
with $\xi'$ as in Lemma \ref{lem:3.7};
\item call $\ol{\la}^{[n]}(\e,\be_0)=
(\ol{\la}^{[n]}_{1}(\e,\be_0),\ldots,\ol{\la}^{[n]}_{r}(\e,\be_0))$
the eigenvalues of $\del_{\be_0}^{2}\ol{L}^{[n]}(\e,\be_0)$;
\item define
\begin{equation} \label{eq:3.10}
\ol{\MM}^{[n]}(x;\e,\be_{0})=
\del_{\be_0}^{2}\ol{L}^{[n]}(\e,\be_0) + \ol{c}_{n}(\e,\be_{0}) + \ol{\gotR}^{[n]}(x;\e,\be_0) .
\end{equation}
\end{enumerate}

We write $\ol{\MM}^{\x}=\{\ol{\MM}^{[n]}(x;\e,\be_{0})\x_{n}(\ol{\la}^{[n]}(\e,\be_{0}))\}_{n\ge-1}$. 
The following result can be proved by reasoning as in the proof of Lemma \ref{lem:3.4}.

%%%%%%%%%%%%%%%%%%%%%%%%%%%%%%%%%%%%%%%%%%%%%%%%%
\begin{lemma}\label{lem:3.15} 
Assume $\ol{\MM}^{\xi}$ to satisfy property 1-$p$. Then for any $0\le n\le p$ one has 
\begin{equation} \nonumber 
\max_{i_{1}+\ldots+i_{r} \le 3}
\left| \del_{\be_{01}}^{i_{1}} \ldots \del_{\be_{0r}}^{i_{r}}
\left( \ol{L}^{[n]}(\e,\be_{0}) - \ol{L}^{[n-1]}(\e,\be_{0}) \right) \right| \le
|\e|^{2} K_{0} \, {\rm e}^{-K_{1} 2^{m_{n}}},
\end{equation}
for suitable $p$-independent constants $K_{0}$ and $K_{1}$.
\end{lemma}
%%%%%%%%%%%%%%%%%%%%%%%%%%%%%%%%%%%%%%%%%%%%%%%%%

Lemma \ref{lem:3.15} implies that, under the same assumptions, one has
\begin{equation} \nonumber
\left\| \ol{\MM}^{[n]} (x;\e,\be_{0}) - \ol{\MM}^{[n-1]} (x;\e,\be_{0}) \right\|
\le |\e|^{2} K_{2} {\rm e}^{-K_{1} 2^{m_{n}}} ,
\end{equation} 
for a suitable $p$-independent constant $K_{2}$.

%%%%%%%%%%%%%%%%%%%%%%%%%%%%%%%%%%%%%%%%%%%%%%%%%
\begin{lemma}\label{lem:3.16} 
$\ol{\MM}^{\x}$ satisfies property 1.
\end{lemma}
%%%%%%%%%%%%%%%%%%%%%%%%%%%%%%%%%%%%%%%%%%%%%%%%%

%%%%%%%%%%%%%%%%%%%%%%%%%%%%%%%%%%%%%%%%%%%%%%%%%
\prova
We shall prove that $\ol{\MM}^{\x}$ satisfies property 1-$p$ for all
$p\ge0$, by induction on $p$. Property 1-0 is trivially satisfied for
$\e$ small enough. Assume $\ol{\MM}^{\x}$ to satisfies
property 1-$p$. Then \eqref{eq:3.9} implies that,
for all $n\le p$ and $x$ such that $\Psi_n(x)\neq0$, one has
\begin{eqnarray}
& & \hskip-.4truecm
\Bigl\| x^2 \uno - \partial_{\be_0}^{2} \ol{L}^{[n-1]}(\e,\be_0) \, \x_{n-1}(\ol{\la}^{[n-1]}(\e,\be_{0}))
- \ol{c}_{n-1}(\e,\be_0) \, \x_{n-1}(\ol{\la}^{[n-1]}(\e,\be_{0})) \nonumber \\
& & \hskip-.4truecm
\qquad\qquad - \, \gotR^{[n-1]}(x,\e,\be_0) \, \x_{n-1}(\ol{\la}^{[n-1]}(\e,\be_{0})) \Bigr\| \nonumber \\
& & \hskip-.4truecm
\quad \ge  \left\| x^2 \uno - \partial_{\be_0}^{2} \ol{L}^{[n-1]}(\e,\be_0) \, 
\x_{n-1}(\ol{\la}^{[n-1]}(\e,\be_{0})) \right\| 
- \, \left\| \ol{c}_{n-1}(\e,\be_{0}) \right\| - \left\| \ol{\gotR}^{[n-1]}(x,\e,\be_0) \right\| \nonumber \\
& & \hskip-.4truecm
\quad \ge x^{2} - \max_{1\le i \le r} \ol{\la}^{[n-1]}_{i}(\e,\be_0) \, 
\x_{n-1}(\ol{\la}^{[n-1]}(\e,\be_{0})) + O({\rm e}^{-\xi'2^{m_{n+1}}}) + O(\e^2 x^{2}) \ge \frac{x^{2}}{2} 
\nonumber
\end{eqnarray}
and hence $\ol{\MM}^{\x}$ satisfies property 1-$(p+1)$.\EP
%%%%%%%%%%%%%%%%%%%%%%%%%%%%%%%%%%%%%%%%%%%%%%%%%

Define
\begin{subequations}\label{eq:3.11}
\begin{align}
& \ol{b}(\al ;\e,\be_{0}):=\lim_{n\to\io} \ol{b}^{[n]}(\al;\e,\be_{0}) ,
\label{eq:3.11a} \\
& \ol{L}^{[\io]}(\e,\be_{0}):=\lim_{n\to\io} \ol{L}^{[n]}(\e,\be_{0}) .
\label{eq:3.11b}
\end{align}
\end{subequations}
By Lemma \ref{lem:3.16}, the limits in (\ref{eq:3.11}) are well defined functions, which are
$C^{\io}$ in both $\e$ and $\be_{0}$ and $2\pi$-periodic in each component of $\be_{0}$.
Moreover $\ol{b}(\al ;\e,\be_{0})$ is analytic in $\al\in\TTT^{d}_{\x''}$ for some $\x''<\x'$.

Set $\ol{G}(\e,\be_{0}):=\del_{\be_{0}}\ol{L}^{[\io]}(\e,\be_{0})$.
The following result is the analogous of Lemma \ref{lem:3.13} 

%%%%%%%%%%%%%%%%%%%%%%%%%%%%%%%%%%%%%%%%%%%%%%%%%
\begin{lemma}\label{lem:3.17} 
For $\e$ small enough there is at least one value
$\ol{\be}_{0}(\e)\in\TTT^{r}$ such that $\ol{G}(\e,\ol{\be}_{0}(\e))=0$ and
$\del_{\be_{0}}^2\ol{L}^{[\io]}(\e,\ol{\be}_{0}(\e)) \le 0$.
\end{lemma} 
%%%%%%%%%%%%%%%%%%%%%%%%%%%%%%%%%%%%%%%%%%%%%%%%%

%%%%%%%%%%%%%%%%%%%%%%%%%%%%%%%%%%%%%%%%%%%%%%%%%
\begin{lemma}\label{lem:3.18} 
Let $\ol{\be}_{0}(\e)$ be as in Lemma \ref{lem:3.17}. Then for $\e$ small enough
and all $n\ge0$ one has $\x_{n}(\la^{[n]}(\e,\be_{0}))\equiv 1$ for all $\be_0$
such that $|\be_{0}-\ol{\be}_{0}(\e)|\le \al_{m_{n+1}}^{2}(\om)$.
\end{lemma} 
%%%%%%%%%%%%%%%%%%%%%%%%%%%%%%%%%%%%%%%%%%%%%%%%%

%%%%%%%%%%%%%%%%%%%%%%%%%%%%%%%%%%%%%%%%%%%%%%%%%
\prova 
Fix ${\be}_{0}=\ol{\be}_{0}(\e)$ so that $\del_{\be_{0}}^2\ol{L}^{[\io]}(\e,\ol{\be}_{0}(\e)) \le 0$
by Lemma \ref{lem:3.17}. One has
\begin{equation} \nonumber 
\del_{\be_{0}}^{2}\ol{L}^{[n]}(\e,\be_{0}) =
\left( \del_{\be_{0}}^{2}\ol{L}^{[n]}(\e,\be_{0}) - 
\del_{\be_{0}}^{2}\ol{L}^{[\io]}(\e,\be_{0}) \right) +
\del_{\be_{0}}^{2}\ol{L}^{[\io]}(\e,\be_{0}) ,
\end{equation}
where, by Lemma \ref{lem:3.15},
\begin{equation} \nonumber 
\left\| \del_{\be_{0}}^{2}\ol{L}^{[n]}(\e,\be_{0}) - 
\del_{\be_{0}}^{2}\ol{L}^{[\io]}(\e,\be_{0}) \right\| \le
K_{3}|\e| {\rm e}^{-K_{1}\,2^{m_{n+1}}} \le \frac{\al_{m_{n+1}}^{2}(\oo)}{2^{14}}, 
\end{equation}
for a suitable positive constant $K_{3}$. Moreover
\begin{equation} \nonumber
\left\| \ol{L}^{[n]}(\e,\be_0) - \ol{L}^{[n]}(\e,\be_{0}(\e)) \right\| \le
K_{4} |\e|^{2} \left| \be_{0}-\ol{\be}_{0}(\e) \right| ,
\end{equation}
for a suitable positive constant $K_{4}$. Therefore,
for $|\be_{0}-\ol{\be}_{0}(\e)|\le \al_{m_{n+1}}^{2}(\om)$, one has
\begin{equation} \nonumber 
\ol{\la}^{[n]}_{j}(\e,\be_{0}) \le 
\ol{\la}^{[\io]}_{j}(\e,\ol{\be}_{0}(\e)) + \frac{\al_{m_{n+1}}^{2}(\oo)}{2^{13}} 
\le \frac{\al_{m_{n+1}}^{2}(\oo)}{2^{13}} , \qquad j=1,\ldots,r,
\end{equation} 
so that the assertion follows. \EP 
%%%%%%%%%%%%%%%%%%%%%%%%%%%%%%%%%%%%%%%%%%%%%%%%%

%%%%%%%%%%%%%%%%%%%%%%%%%%%%%%%%%%%%%%%%%%%%%%%%%
\begin{lemma}\label{lem:3.19} 
One can choose the matrices $\ol{c}_{n}(\e,\be_{0})$ and $\ol{\gotR}^{[n]}(x;\e,\be_{0})$ 
in \eqref{eq:3.10} so that
the function $\be(t)=\ol{\be}_{0}(\e)+\ol{b}(\al;\e,\ol{\be}_{0}(\e))$ solves \eqref{eq:1.5}. 
\end{lemma} 
%%%%%%%%%%%%%%%%%%%%%%%%%%%%%%%%%%%%%%%%%%%%%%%%%

%%%%%%%%%%%%%%%%%%%%%%%%%%%%%%%%%%%%%%%%%%%%%%%%%
\prova
We want to prove by induction that we can choose the matrices
$\ol{c}_{n}(\e,\be_{0})$ and $\ol{\gotR}^{[n]}(x;\e,\be_{0})$, for $n\ge 0$, so as to have
\begin{equation} \label{eq:3.12}
\ol{\MM}^{[n]}(x;\e,\ol{\be}_0(\e))=\MM^{[n]}(x;\e,\ol{\be}_0(\e))
\end{equation}
for all $n\ge -1$. This will imply that, when setting $\be_0=\ol{\be}_{0}(\e)$,
$\MM$ satisfies property 1 as well as $\ol{\MM}$. Hence we obtain
$b(\al;\e,\ol{\be}_{0}(\e))=\ol{b}(\al;\e,\ol{\be}_{0}(\e))$ and we can apply
Lemma \ref{lem:3.14} with $\be_{0}(\e)=\ol{\be}_0(\e)$ to deduce that
$\ol{\be}_{0}(\e)+ \ol{b}(\al;\e,\ol{\be}_{0}(\e))$ solves the range equation \eqref{eq:1.7a}.

For $n=-1$ the identity \eqref{eq:3.12} trivially holds.
Assume that \eqref{eq:3.12} is satisfied up to $n<p$ for a suitable choice of the matrices
$\ol{c}_{n}(\e,\be_{0})$ and $\ol{\gotR}^{[n]}(x;\e,\be_{0})$, with $0 \le n < p$.
Then, for $|\be_0-\ol{\be}_{0}(\e)| \le \al_{m_{n+1}}^{2}(\om)$,
$\MM$ satisfies property 1-$p$ and hence $\MM^{[p]}(x;\e,\be_0)$ is well defined.
Moreover for such $\be_0$ one has $\ol{L}^{[p]}(\e,\be_{0})=L^{\le p}(\e,\be_{0})$ and hence
\begin{equation} \nonumber
\del_{\be_0}^{2} \ol{L}^{[p]}(\e,\be_{0}) = \MM^{[p]}(0;\e,\be_0) - d_{p}(\e,\be_0) ,
\end{equation}
with the matrix $d_{p}(\e,\be_0)=\del_{\be_0}c_p'(\e,\be_0)+c_p(\e,\be_0)$ such that
$\|d_{p}(\e,\be_{0})\| \le K_{5}|\e|^{2} {\rm e}^{-\x' 2^{m_{p+1}}}$,
for a suitable positive constant $K_{5}$. Then \eqref{eq:3.10} gives
\begin{equation} \nonumber
\ol{\MM}^{[p]}(x;\e,\be_0) = \MM^{[p]}(0;\e,\be_0) - d_{p}(\e,\be_0) +
\ol{c}_{p}(\e,\be_0) + \ol{\gotR}^{[p]}(x;\e,\be_0) ,
\end{equation}
and one can take the matrices $\ol{c}_{p}(\e,\be_0)$ and $\ol{\gotR}^{[p]}(x;\e,\be_0)$ so that
\begin{equation} \nonumber
\ol{c}_{p}(\e,\be_0) = d_{p}(\e,\be_0) , \qquad
\ol{\gotR}^{[p]}(x;\e,\be_0) = \MM^{[p]}(x;\e,\be_0) - \MM^{[p]}(0;\e,\be_0) ,
\end{equation}
for $|\be_0-\ol{\be}_{0}(\e)| \le \al_{m_{n+1}}^{2}(\om)$. Thanks to Lemma \ref{lem:3.5},
this yields that both inequalities \eqref{eq:3.9} are satisfied.
Therefore \eqref{eq:3.12} follows for $n=p$.

On the other hand one has
\begin{equation} \nonumber
G(\e,\ol{\be}_0(\e))= \lim_{n\to\io}\del_{\be_{0}}L^{\le n}(\e,\ol{\be}_{0}(\e)) =
\del_{\be_0} \ol{L}^{[\io]}(\e,\ol{\be}_{0}(\e)) = \ol{G}(\e,\ol{\be}_{0}(\e)) = 0
\end{equation}
and hence the bifurcation equation \eqref{eq:1.7b} is solved as well.
\EP
%%%%%%%%%%%%%%%%%%%%%%%%%%%%%%%%%%%%%%%%%%%%%%%%%

%%%%%%%%%%%%%%%%%%%%%%%%%%%%%%%%%%%%%%%%%%%%%%%%%
%%%%%%%%%%%%%%%%%%%%%%%%%%%%%%%%%%%%%%%%%%%%%%%%%
% References 
%%%%%%%%%%%%%%%%%%%%%%%%%%%%%%%%%%%%%%%%%%%%%%%%%
%%%%%%%%%%%%%%%%%%%%%%%%%%%%%%%%%%%%%%%%%%%%%%%%%

\end{document}